 
 
 
 

  \ifx\osumriploaded\MYundefined  
  \else
    \immediate\write16{osumrip.sty ALREADY loaded.}
  \endinput\fi

  \def\osumriploaded{\relax}

  \magnification=1200 
  \hoffset=.35truein  
  \hsize=355pt 
  \vsize = 500 pt
  \baselineskip=11pt 
  \lineskip=1.1pt
  \lineskiplimit=.8pt
  \topskip=12pt 
  \bigskipamount=10pt plus 4pt minus 1pt
  \hfuzz=\hsize
  \overfullrule = 0 pt 
  \vbadness=10000
  \hbadness=10000
  \nopagenumbers
  \widowpenalty=5000
  \pageno=1
  \footline{\ifnum \pageno>0 \hss\tenrm\folio\hss\fi}
  \mathsurround=1pt
  \def\prose{\kern\mathsurround}

  \def\StdPretolerance{100}
  \tolerance=\StdPretolerance

  \def\StdParskip{0pt}   
  \parskip=\StdParskip
  \parindent=0.5cm
 
  \def\tenpoint{}
 
 
  \font\Bigbf=cmbx10 at 14.4pt 
  \font\twelvebf=cmbx12
  \font\tenbf=cmbx10
  
  \font\tensmc=cmcsc10

 \def\smc{\tensmc}
 
 \def \Smallfonts {\relax } 
  \ifx\amspptloaded@AmS\relax
     \def\Smallfonts{\eightpoint}
  \fi
  \def\next{AMSPPT}
  \ifx\styname\next 
     \def\Smallfonts{\eightpoint}
  \fi

  \def\Hfont{\Bigbf}
  \def\Authorfont{\twelvebf}
  \def\HHfont{\twelvebf} 
  \def\HHHfont{\tenbf}
  \def\Figurefont{\bf}
  \def\Bibfont{\tenbf}
  
 \def \thfont {\smc }
 \def \pffont {\smc }
 \def \rkfont {\smc }
 \def \dffont {\smc }
 \def \egfont {\smc }
 

 
  \chardef\CatAt\the\catcode`\@
  \catcode`\@=11

 \let\wlog@ld\wlog 
 \def\wlog#1{\relax}

   \def\hexnumber@#1{\ifcase#1 
     0\or1\or2\or3\or4\or5\or6\or7\or8\or9\or 
     A\or B\or C\or D\or E\or F\fi}

 
 \newif\ifIN@
 \def\m@rker{\m@@rker}
 \def\IN@{\expandafter\INN@\expandafter}
 \long\def\INN@0#1@#2@{\long\def\NI@##1#1##2##3\ENDNI@
    {\ifx\m@rker##2\IN@false\else\IN@true\fi}%
     \expandafter\NI@#2@@#1\m@rker\ENDNI@}

  \newtoks\Initialtoks@  \newtoks\Terminaltoks@
  \def\SPLIT@{\expandafter\SPLITT@\expandafter}
  \def\SPLITT@0#1@#2@{\def\TTILPS@##1#1##2@{%
     \Initialtoks@{##1}\Terminaltoks@{##2}}\expandafter\TTILPS@#2@}


  \newtoks\Trimtoks@

 \def\ForeTrim@{\expandafter\ForeTrim@@\expandafter}
 \def\ForePrim@0 #1@{\Trimtoks@{#1}}
 \def\ForeTrim@@0#1@{\IN@0\m@rker. @\m@rker.#1@%
     \ifIN@\ForePrim@0#1@%
     \else\Trimtoks@\expandafter{#1}\fi}

  \def\Trim@0#1@{%
      \ForeTrim@0#1@%
      \IN@0 @\the\Trimtoks@ @%
        \ifIN@ 
             \SPLIT@0 @\the\Trimtoks@ @\Trimtoks@\Initialtoks@
             \IN@0\the\Terminaltoks@ @ @%
                 \ifIN@
                 \else \Trimtoks@ {FigNameWithSpace}%
                 \fi
        \fi
      }

 
 
 \def\Hrule{\hrule width0pt height0pt}

 \newskip\LastSkip
 \def\SaveLastSkip{\LastSkip\lastskip}

 \def\NoindentAfter{\everypar={\setbox0=\lastbox\everypar={}}}

 \long\def\H#1\par#2\par{\notenumber=0%
    \hbox to 0pt{}\vskip15pt plus 1pt 
    {\baselineskip=15pt plus 1pt\parindent=0pt\parskip=0pt\frenchspacing
    \leftskip=0pt plus .2\hsize minus .1\hsize 
    \rightskip=0pt plus .2\hsize minus .1\hsize 
    \def\\{\unskip\break} 
    \pretolerance=10000 \Hfont #1\unskip\break
    \vskip20pt plus 1pt\Authorfont #2\unskip\break\par}
    \vskip20pt plus 1pt%
    \NoindentAfter\par\rm
    }




 \newdimen\PageRemainder
  \def\SetPageRemainder{\PageRemainder=\pagegoal\advance\PageRemainder by
   -1\pagetotal}
 
  \def\Rpt@{}\def\Rpt@@{}
  
  \long\def\HH#1\par{\par
  \SaveLastSkip\removelastskip\goodbreak
  \ifdim\LastSkip<24pt
     \LastSkip 24pt plus 1pt\fi
  \SetPageRemainder\advance\PageRemainder-\LastSkip
  \ifdim\PageRemainder<150pt
       \edef\Rpt@{remain = \the\PageRemainder\noexpand\\
                pagetotal=\the\pagetotal\noexpand\\
                           pagegoal=\the\pagegoal}%
          \else\def\Rpt@{}\fi
   \ifdim\PageRemainder<72pt 
             \edef\Rpt@@{\noexpand\\
                      Had HH PageRemainder$<$72pt\noexpand\\
                      Hence forced break!}%
     \vskip 0pt plus .1\PageRemainder\eject 
    \else\edef\Rpt@@{}
    \fi
    \vskip\LastSkip\Hrule 
    \pretolerance=10000\rightskip=0pt plus 3em
    \hangafter1 \hangindent=2.2em%
    \noindent
    \HHfont \unskip \Ednote{\Rpt@\Rpt@@}%
            \def\Rpt@{}\def\Rpt@@{}%
            \ignorespaces
            #1\par\rightskip=0pt\pretolerance=\StdPretolerance%
      \NoindentAfter\tenpoint\rm\vskip 8pt plus 1pt minus 1pt}
 
  \long\def\HHH#1\par{\par%
  \SaveLastSkip\removelastskip\goodbreak
  \ifdim\LastSkip<10pt
     \LastSkip 10pt plus 1pt\fi
  \SetPageRemainder\advance\PageRemainder-\LastSkip
  \ifdim\PageRemainder<150pt
       \edef\Rpt@{remain = \the\PageRemainder\noexpand\\
                pagetotal=\the\pagetotal\noexpand\\
                           pagegoal=\the\pagegoal}%
          \else\def\Rpt@{}\fi
   \ifdim\PageRemainder<60pt  
             \edef\Rpt@@{\noexpand\\
                      Had HHH PageRemainder$<$60pt\noexpand\\
                      Hence forced break!}%
       \vskip 0pt plus .1\PageRemainder\eject 
   \else\edef\Rpt@@{}
   \fi
   \vskip\LastSkip\Hrule\par\noindent
   \HHHfont \unskip\Ednote{\Rpt@\Rpt@@}\ignorespaces
   #1\unskip.\quad\rm\ignorespaces
   \ignorepars}
       
  \long\def\ignorepars#1\par{\def\Test{#1}%
     \ifx\Test\Empty\def\This{\ignorepars}%
        \else\def\This{\Test\par}\fi
           \This}

 
 \def\ProcBreak{\par\ifdim\lastskip<8pt
    \removelastskip
    \penalty-200\vskip8pt plus1pt\fi}

 \def \th #1\par{\ProcBreak \noindent
   {\thfont\ignorespaces #1\unskip.}\quad\it}
 \def \endth {\ProcBreak\rm }
 
 
 \def\pf #1\par{\ProcBreak %
    \noindent\pffont#1\unskip:\kern.5pt
       \raise-.4pt\hbox{---}\quad\rm}
 

  \def\qedbox{\hbox{\vbox{
    \hrule width0.2cm height0.2pt 
    \hbox to 0.2cm{\vrule height 0.2cm width 0.2pt 
             \hfil\vrule height0.2cm width 0.2pt}
    \hrule width0.2cm height 0.2pt}}}

  \def\qed{\ifmmode\qedbox
    \else\hglue4mm\unskip\hfill\qedbox\ProcBreak\fi} 

  \def \rk #1\par{\ProcBreak
     \noindent{\rkfont\ignorespaces #1\unskip.}\quad\rm}
  
  \def \endrk {\ProcBreak }
 
  \def \df #1\par{\ProcBreak 
     \noindent{\dffont\ignorespaces #1\unskip.}\quad\rm}

  \def \eg #1\par{\ProcBreak 
     \noindent{\egfont\ignorespaces #1\unskip.}\quad\rm}



  \newdimen\Overhang

   \def\MaxTag@#1#2#3#4#5{\setbox0=\hbox{#4\ignorespaces#2\unskip}%
     \dimen0=\wd0\advance\dimen0 by#3
     \ifdim\dimen0<#5\relax\dimen0=#5\fi
     \expandafter\edef\csname #1Hang\endcsname{\the\dimen0}}

 \def\MaxItemTag#1{\MaxTag@{Item}{#1}{.4em}{\ItemStyle}{\parindent}}%
 \def\MaxItemItemTag#1{%
        \MaxTag@{ItemItem}{#1}{.4em}{\ItemItemStyle}{\parindent}}
 \def\MaxNrTag#1{\MaxTag@{Nr}{#1}{.5em}{\NrStyle}{\parindent}}
 \def\MaxReferenceTag#1{%
        \MaxTag@{Reference}{[#1]}{.6em}{\Smallfonts}{\parindent}}
 \def\MaxFootTag#1{\MaxTag@{Foot}{#1}{.4em}{\Smallfonts}{\z@}}

  \def\SetOverhang@{\Overhang=.8\dimen0%
     \advance\Overhang by \wd0\relax
     \ifdim\Overhang>\hangindent\relax
       \advance\Overhang by .25\dimen0%
       \Ednote{Tag is pushing text.}
     \else\Overhang=\hangindent
     \fi}

   \def\Item#1{\par\noindent
      \hangafter1\hangindent=\ItemHang
      \setbox0=\hbox{\ItemStyle\ignorespaces#1\unskip}%
      \dimen0=.4em\SetOverhang@
      \rlap{\box0}\kern\Overhang\ignorespaces}

   \def\ItemStyle{\rm}
   \MaxItemTag{(iii)}
         
   \def\ItemItem#1{\par\noindent
      \hangafter1\hangindent=\ItemItemHang
      \setbox0=\hbox{\ItemItemStyle\ignorespaces#1\unskip}%
      \dimen0=.4em\SetOverhang@
      \advance\hangindent by \ItemHang
      \kern\ItemHang\rlap{\box0}%
      \kern\Overhang\ignorespaces}

    \def\ItemItemStyle{\rm}
    \MaxItemItemTag{(iii)}

  \def\Nr#1{\par\noindent\hangindent=\NrHang 
    \setbox0=\hbox{\NrStyle\ignorespaces#1\unskip}%
    \dimen0=.5em\SetOverhang@
    \rlap{\box0}\kern\Overhang
    \hangindent=\z@\ignorespaces}

   \def\NrStyle{\rm}
   \MaxNrTag{(2)}

   \newskip\Rosterskip\Rosterskip 1pt plus1pt 
   \def\Roster{\par\ifdim\lastskip<\Rosterskip\removelastskip\vskip\Rosterskip\fi
    \bgroup}
   \def\endRoster{\par\global\edef\LastSkip@{\the\lastskip}\removelastskip
       \egroup\penalty-50\LastSkip\LastSkip@\relax
       \ifdim\LastSkip<\Rosterskip\LastSkip\Rosterskip\fi
       \vskip\LastSkip}




\def\cite#1{
    \def\nextiii@##1,##2\end@{{\frenchspacing\bf
      \lBr\ignorespaces##1\unskip{\rm,~\ignorespaces##2}\rBr}}%
    \IN@0,@#1@%
    \ifIN@\def\next{\nextiii@#1\end@}\else
    \def\next{{\bf\lBr#1\rBr}}\fi\next}



  
   \def \Bib#1\par{%
       \par\removelastskip\SetPageRemainder
       \ifdim\PageRemainder < 97pt
        \ifdim\PageRemainder > 0pt
        \vfill\eject
       \fi\fi
    \ProcBreak \par\begingroup\parskip=0 pt%
    \goodbreak \vskip 15 pt plus 10 pt
    \noindent\null\hfill\Bibfont
      \ignorespaces #1\unskip\hfill\null\par 
    \frenchspacing \Smallfonts\rm
    \parskip=2.5 pt plus 1 pt minus.5pt%
    \nobreak\vskip 12pt plus 2pt minus2pt\nobreak
    \leftskip=0 pt \baselineskip=10.5pt}


 \def\ReferenceTagSlide{0em}
  \def\ReferenceTagGap{.5em}
  \def\ReferenceHang{30pt}

  \def \rf#1{\par\noindent
     \hangafter1\hangindent=\ReferenceHang      
     \setbox0=\hbox{\Smallfonts[\ignorespaces#1\unskip]}%
     \dimen0=\ReferenceTagGap\SetOverhang@
     \rlap{\kern\ReferenceTagSlide\box0}%
     \kern\Overhang\ignorespaces}

  \def\ref#1\par#2\par#3\par#4\par{%
     \rf{#1}#2\unskip,\ #3\unskip,\
     #4\unskip.}

  \def\endBib{\par\endgroup\vskip 12pt minus 6pt }


  \def\Coordinates{\bigskip\bigskip
     \vtop\bgroup\leftskip=\parindent
    \parskip=3pt \parindent=0pt\pretolerance=10000%
    \def\\{\hfil \break} \frenchspacing\rm }
    
  \def\endCoordinates{\par\egroup}


  \newcount\notenumber
  
  \def\note{\advance\notenumber by 1
    \footnote{\the\notenumber)}}

  \def\footnote#1{\let\@sf\empty
    \ifhmode\edef\@sf{\spacefactor\the\spacefactor}\/\fi
    \sam${}^{\fam0 #1}$\@sf\vfootnote{#1}}%

  \def\vfootnote#1{\insert\footins\bgroup
     \interlinepenalty100 \splittopskip=1pt
     \floatingpenalty=20000
     \leftskip=0pt\rightskip=0pt%
     \parindent=.3em
     \Smallfonts\rm
     \FootItem@{#1}
     \futurelet\next\fo@t}

  \def\FootItem@#1{\par\hangafter1\hangindent=\FootHang
     \setbox0=\hbox{\ignorespaces#1\unskip}%
     \dimen0=.4em\SetOverhang@
     \noindent\rlap{\box0}\kern\Overhang\ignorespaces}

  \MaxFootTag{2)}

  \def\fo@t{\ifcat\bgroup\noexpand\next \let\next\f@@t
    \else\let\next\f@t\fi \next}
  \def\f@@t{\bgroup\aftergroup\@foot\let\next}
  \def\f@t#1{\baselineskip=10pt\lineskip=1pt
            \lineskiplimit=0pt #1\@foot}%
  \def\@foot{
        \hbox{\vrule height0pt depth5pt width0pt}
        \egroup}
  \skip\footins=12 pt plus 0pt minus 0pt 
  \count\footins=1000 
  \dimen\footins=8in 


  \let\Ninepoint\ninepoint

  \def\prose{\kern1.5pt}

 \def \Blackbox
   {\leavevmode\hskip .3pt \vbox
   {\hrule height 5pt\hbox{\hskip 4.5pt}}\hskip .5pt}

 \def \XX{\Blackbox\kern.5pt\Blackbox} 

  \def\.{.\kern1pt}

    \def\Hyphen{\edef\this{\the\hyphenchar\font}%
          \hyphenchar\font=-1\char\this\hyphenchar\font=\this}

 \ifx\undefined\text
  \def\text#1{\hbox{\rm #1}}\fi 


 
 \newcount\Ht 

 \def \Acc{\expandafter }

 \def\swthat{\raise -1.1 ex\hbox{$\widehat{}$}}
 \def\swttilde{\raise -1.2 ex\hbox{$\widetilde{}$}}
 \def \overdot{{\raise .2 ex \hbox to 0pt {\hss\bf\smash{.}\hss}}}
 \def \overcircle{{\raise .1 ex \hbox to 0pt
    {\mathsurround=0pt$\eightpoint\scriptstyle\hss\circ\hss$}}}

 \def \Mathaccent#1#2{{\mathsurround=0 pt
  \setbox4=\hbox{$\vphantom{#2}$}
  \Ht=\ht4 
  \setbox5=\hbox{${#1}$}
  \setbox6=\hbox{${#2}$}
  \setbox7=\hbox to .5\wd6{}
  \copy7\kern .1\Ht \raise\Ht sp\hbox{\copy5}\kern-.1\Ht 
  \copy7\llap{\box6}
  }}

  \def\SwtCheck #1{
  \ifmmode \check{#1}%
    \else \v {#1}%
    \fi}

 \def\barpartial {%
   \kern .17 em 
    \overline {\kern -.17 em\partial\kern-.03 em}%
    \kern .03 em}

 \def\Overline#1{\setbox1=\hbox{\sam ${#1}$}%
      \ifdim \wd1 > 6pt
    \kern .11 em
    \overline {\kern -.11 em#1\kern-.14 em}
    \kern .14 em 
  \else
    \kern .03 em 
    \overline {\kern -.03 em#1\kern-.04 em}
    \kern .04 em 
  \fi}

 \def\SOverline#1{\setbox1=\hbox{\sam ${#1}$}%
      \ifdim \wd1 > 7pt
    \kern .22 em 
    \overline {\kern -.22 em#1\kern-.09 em}%
    \kern .09 em 
  \else
    \kern .10 em 
    \overline {\kern -.10 em#1\kern-.04 em}%
    \kern .04 em 
  \fi}


 \def\Underline#1{\setbox1=\hbox{\sam ${#1}$}%
      \ifdim \wd1 > 6pt
    \kern .11 em
    \underline {\kern -.11 em#1\kern-.14 em}
    \kern .14 em 
  \else
    \kern .03 em 
    \underline {\kern -.03 em#1\kern-.04 em}
    \kern .04 em 
  \fi}

 \def\SUnderline#1{\setbox1=\hbox{\sam ${#1}$}%
      \ifdim \wd1 > 7pt
    \kern .04 em 
    \underline {\kern -.04 em#1\kern-.2 em}%
    \kern .2 em 
  \else
    \kern .0 em 
    \underline {\kern -.0 em#1\kern-.15 em}%
    \kern .15 em 
  \fi}


 \ifx\MYUNDEFINED\BoxedEPSF
   \let\temp\relax
 \else
   \message{}
   \message{ !!! BoxedEPS %
         or BoxedArt macros already defined !!!}
   \let\temp\endinput
 \fi
  \temp
 
 \chardef\CatAt\the\catcode`\@
 \catcode`\@=11
 \chardef\C@tColon\the\catcode`\:
 \chardef\C@tSemicolon\the\catcode`\;
 \chardef\C@tQmark\the\catcode`\?
 \chardef\C@tEmark\the\catcode`\!

 \def\PunctOther@{\catcode`\:=12
   \catcode`\;=12 \catcode`\?=12 \catcode`\!=12}
 \PunctOther@

 \let\wlog@ld\wlog 
 \def\wlog#1{\relax} 

 \newif\ifIN@
 \newdimen\XShift@ \newdimen\YShift@ 
 \newtoks\Realtoks
 
  %
 \newdimen\Wd@ \newdimen\Ht@
 \newdimen\Wd@@ \newdimen\Ht@@
 \newdimen\TT@
 \newdimen\LT@
 \newdimen\BT@
 \newdimen\RT@
 \newdimen\XSlide@ \newdimen\YSlide@ 
 \newdimen\TheScale  
 \newdimen\FigScale  
 \newdimen\ForcedDim@@

 \newtoks\EPSFDirectorytoks@
 \newtoks\EPSFNametoks@
 \newtoks\BdBoxtoks@
 \newtoks\LLXtoks@  
 \newtoks\LLYtoks@

 \newif\ifNotIn@
 \newif\ifForcedDim@
 \newif\ifForceOn@
 \newif\ifForcedHeight@
 \newif\ifPSOrigin

 \newread\EPSFile@ 
 
  \def\ms@g{\immediate\write16}

 \newif\ifIN@\def\IN@{\expandafter\INN@\expandafter}
  \long\def\INN@0#1@#2@{\long\def\NI@##1#1##2##3\ENDNI@
    {\ifx\m@rker##2\IN@false\else\IN@true\fi}%
     \expandafter\NI@#2@@#1\m@rker\ENDNI@}
  \def\m@rker{\m@@rker}

  \newtoks\Initialtoks@  \newtoks\Terminaltoks@
  \def\SPLIT@{\expandafter\SPLITT@\expandafter}
  \def\SPLITT@0#1@#2@{\def\TTILPS@##1#1##2@{%
     \Initialtoks@{##1}\Terminaltoks@{##2}}\expandafter\TTILPS@#2@}


  \newtoks\Trimtoks@

 \def\ForeTrim@{\expandafter\ForeTrim@@\expandafter}
 \def\ForePrim@0 #1@{\Trimtoks@{#1}}
 \def\ForeTrim@@0#1@{\IN@0\m@rker. @\m@rker.#1@%
     \ifIN@\ForePrim@0#1@%
     \else\Trimtoks@\expandafter{#1}\fi}

  \def\Trim@0#1@{%
      \ForeTrim@0#1@%
      \IN@0 @\the\Trimtoks@ @%
        \ifIN@ 
             \SPLIT@0 @\the\Trimtoks@ @\Trimtoks@\Initialtoks@
             \IN@0\the\Terminaltoks@ @ @%
                 \ifIN@
                 \else \Trimtoks@ {FigNameWithSpace}%
                 \fi
        \fi
      }


   \newtoks\pt@ks
   \def \getpt@ks 0.0#1@{\pt@ks{#1}}
   \dimen0=0pt\relax\expandafter\getpt@ks\the\dimen0@

  \newtoks\Realtoks
  \def\Real#1{%
    \dimen2=#1%
      \SPLIT@0\the\pt@ks @\the\dimen2@
       \Realtoks=\Initialtoks@
            }

   \newdimen\Product
   \def\Mult#1#2{%
     \dimen4=#1\relax
     \dimen6=#2%
     \Real{\dimen4}%
     \Product=\the\Realtoks\dimen6%
        }

 \newdimen\Inverse
 \newdimen\hmxdim@ \hmxdim@=8192pt
 \def\Invert#1{%
  \Inverse=\hmxdim@
  \dimen0=#1%
  \divide\Inverse \dimen0%
  \multiply\Inverse 8}

   \def\Rescale#1#2#3{
              \divide #1 by 100\relax
              \dimen2=#3\divide\dimen2 by 100 \Invert{\dimen2}%
              \Mult{#1}{#2}%
              \Mult\Product\Inverse 
              #1=\Product}

  \def\Scale#1{\dimen0=\TheScale %
      \divide #1 by  1280 
      \divide \dimen0 by 5120 %
      \multiply#1 by \dimen0 
      \divide#1 by 10   
     }
 

 \newbox\scrunchbox

 \def\Scrunched#1{{\setbox\scrunchbox\hbox{#1}%
   \wd\scrunchbox=0pt
   \ht\scrunchbox=0pt
   \dp\scrunchbox=0pt
   \box\scrunchbox}}

 \def\Shifted@#1{%
   \vbox {\kern-\YShift@
       \hbox {\kern\XShift@\hbox{#1}\kern-\XShift@}%
           \kern\YShift@}}


 \def\cBoxedEPSF#1{{{}\leavevmode 
   \ReadNameAndScale@{#1}%
   \SetEPSFSpec@
   \ReadEPSFile@ \ReadBdB@x  
     \TrimFigDims@ 
     \CalculateFigScale@  
     \ScaleFigDims@
     \SetInkShift@
   \hbox{$\mathsurround=0pt\relax
         \vcenter{\hbox{%
             \FrameSpider{\hskip-.4pt\vrule}%
             \vbox to \Ht@{\offinterlineskip\parindent=\z@%
                \FrameSpider{\vskip-.4pt\hrule}\vfil 
                \hbox to \Wd@{\hfil}%
                \vfil
                \InkShift@{\EPSFSpecial{\EPSFSpec@}{\FigSc@leReal}}%
             \FrameSpider{\hrule\vskip-.4pt}}%
         \FrameSpider{\vrule\hskip-.4pt}}}%
     $}%
    \CleanRegisters@ 
    \ms@g{ *** Box composed for the %
         EPSF file \the\EPSFNametoks@}%
    }}      

 \def\tBoxedEPSF#1{\setbox4\hbox{\cBoxedEPSF{#1}}%
     \setbox4\hbox{\raise -\ht4 \hbox{\box4}}%
     \box4
      }

 \def\bBoxedEPSF#1{\setbox4\hbox{\cBoxedEPSF{#1}}%
     \setbox4\hbox{\raise \dp4 \hbox{\box4}}%
     \box4
      }

  \let\BoxedEPSF\cBoxedEPSF

   %

   %
  \def\gLinefigure[#1scaled#2]_#3{%
        \BoxedEPSF{#3 scaled #2}}
    
   %
   \let\EPSFfile\bBoxedEPSF
  
  \def\EPSFxsize{\afterassignment\ForceW@\ForcedDim@@}
      \def\ForceW@{\ForcedDim@true\ForcedHeight@false}
  
  \def\EPSFysize{\afterassignment\ForceH@\ForcedDim@@}
      \def\ForceH@{\ForcedDim@true\ForcedHeight@true}

  %
 \def\ReadNameAndScale@#1{\IN@0 scaled@#1@
   \ifIN@\ReadNameAndScale@@0#1@%
   \else \ReadNameAndScale@@0#1 scaled\DefaultMilScale @
   \fi}
  
 \def\ReadNameAndScale@@0#1scaled#2@{
    \let\OldBackslash@\\%
    \def\\{\OtherB@ckslash}%
    \edef\temp@{#1}%
    \Trim@0\temp@ @%
    \EPSFNametoks@\expandafter{\the\Trimtoks@ }%
    \FigScale=#2 pt%
    \let\\\OldBackslash@
    }
 
 \def\SetDefaultEPSFScale#1{%
      \global\def\DefaultMilScale{#1}}

 \SetDefaultEPSFScale{1000}

  %
 \def \SetBogusBbox@{%
     \global\BdBoxtoks@{ BoundingBox:0 0 100 100 }%
     \global\def\BdBoxLine@{ BoundingBox:0 0 100 100 }%
     \ms@g{ !!! Will use placeholder !!!}%
     }

 {\catcode`\%=12\gdef\P@S@{

 \def\ReadEPSFile@{
     \openin\EPSFile@\EPSFSpec@
     \relax  
  \ifeof\EPSFile@
     \ms@g{}%
     \ms@g{ !!! EPS FILE \the\EPSFDirectorytoks@
       \the\EPSFNametoks@\ WAS NOT FOUND !!!}
     \SetBogusBbox@
  \else
   \begingroup
   \catcode`\%=12\catcode`\:=12\catcode`\!=12
   \catcode`\G=14\catcode`\\=14\relax
   \global\read\EPSFile@ to \BdBoxLine@
   \IN@0\P@S@ @\BdBoxLine@ @%
   \ifIN@ 
     \NotIn@true
     \loop   
       \ifeof\EPSFile@\NotIn@false 
         \ms@g{}%
         \ms@g{ !!! BoundingBox NOT FOUND IN %
            \the\EPSFDirectorytoks@\the\EPSFNametoks@\ !!! }%
         \SetBogusBbox@
       \else\global\read\EPSFile@ to \BdBoxLine@
       \fi
       \global\BdBoxtoks@\expandafter{\BdBoxLine@}%
       \IN@0BoundingBox:@\the\BdBoxtoks@ @%
       \ifIN@\NotIn@false\fi%
     \ifNotIn@\repeat
   \else
         \ms@g{}%
         \ms@g{ !!! \the\EPSFNametoks@\ not PS!\  !!!}%
         \SetBogusBbox@
   \fi
  \endgroup\relax
  \fi
  \closein\EPSFile@ 
   }

  \def\ReadBdB@x{
   \expandafter\ReadBdB@x@\the\BdBoxtoks@ @}
  
  \def\ReadBdB@x@#1BoundingBox:#2@{
    \ForeTrim@0#2@%
    \IN@0atend@\the\Trimtoks@ @%
       \ifIN@\Trimtoks@={0 0 100 100 }%
         \ms@g{}%
         \ms@g{ !!! BoundingBox not found in %
         \the\EPSFDirectorytoks@\the\EPSFNametoks@\space !!!}%
         \ms@g{ !!! It must not be at end of EPSF !!!}%
         \ms@g{ !!! Will use placeholder !!!}%
       \fi
    \expandafter\ReadBdB@x@@\the\Trimtoks@ @%
   }
    
  \def\ReadBdB@x@@#1 #2 #3 #4@{
      \Wd@=#3bp\advance\Wd@ by -#1bp%
      \Ht@=#4bp\advance\Ht@ by-#2bp%
       \Wd@@=\Wd@ \Ht@@=\Ht@ 
       \LLXtoks@={#1}\LLYtoks@={#2}
      \ifPSOrigin\XShift@=-#1bp\YShift@=-#2bp\fi 
     }

   %
   \def\G@bbl@#1{}
   \bgroup
     \global\edef\OtherB@ckslash{\expandafter\G@bbl@\string\\}
   \egroup

  \def\SetEPSFDirectory{
           \bgroup\PunctOther@\relax
           \let\\\OtherB@ckslash
           \SetEPSFDirectory@}

 \def\SetEPSFDirectory@#1{
    \edef\temp@{#1}%
    \Trim@0\temp@ @
    \global\toks1\expandafter{\the\Trimtoks@ }\relax
    \egroup
    \EPSFDirectorytoks@=\toks1
    }

 \def\SetEPSFSpec@{%
     \bgroup
     \let\\=\OtherB@ckslash
     \global\edef\EPSFSpec@{%
        \the\EPSFDirectorytoks@\the\EPSFNametoks@}%
     \global\edef\EPSFSpec@{\EPSFSpec@}%
     \egroup}

  %
 \def\TrimTop#1{\advance\TT@ by #1}
 \def\TrimLeft#1{\advance\LT@ by #1}
 \def\TrimBottom#1{\advance\BT@ by #1}
 \def\TrimRight#1{\advance\RT@ by #1}

 \def\TrimBoundingBox#1{%
   \TrimTop{#1}%
   \TrimLeft{#1}%
   \TrimBottom{#1}%
   \TrimRight{#1}%
       }

 \def\TrimFigDims@{%
    \advance\Wd@ by -\LT@ 
    \advance\Wd@ by -\RT@ \RT@=\z@
    \advance\Ht@ by -\TT@ \TT@=\z@
    \advance\Ht@ by -\BT@ 
    }

  %
  \def\ForceWidth#1{\ForcedDim@true
       \ForcedDim@@#1\ForcedHeight@false}
  
  \def\ForceHeight#1{\ForcedDim@true
       \ForcedDim@@=#1\ForcedHeight@true}

  \def\ForceOn{\ForceOn@true}
  \def\ForceOff{\ForceOn@false\ForcedDim@false}
  
  \def\epsfxsize{\afterassignment\ForceW@\ForcedDim@@}
      \def\ForceW@{\ForcedDim@true\ForcedHeight@false}
  
  \def\epsfysize{\afterassignment\ForceH@\ForcedDim@@}
      \def\ForceH@{\ForcedDim@true\ForcedHeight@true}
  
  \def\CalculateFigScale@{%
     \ifForcedDim@\FigScale=1000pt
           \ifForcedHeight@
                \Rescale\FigScale\ForcedDim@@\Ht@
           \else
                \Rescale\FigScale\ForcedDim@@\Wd@
           \fi
     \fi
     \Real{\FigScale}%
     \edef\FigSc@leReal{\the\Realtoks}%
     }
   
  \def\ScaleFigDims@{\TheScale=\FigScale
      \ifForcedDim@
           \ifForcedHeight@ \Ht@=\ForcedDim@@  \Scale\Wd@
           \else \Wd@=\ForcedDim@@ \Scale\Ht@
           \fi
      \else \Scale\Wd@\Scale\Ht@        
      \fi
      \ifForceOn@\relax\else\global\ForcedDim@false\fi
      \Scale\LT@\Scale\BT@  
      \Scale\XShift@\Scale\YShift@
      }
      
 \def\HideReservedBoxes{\global\def\FrameSpider##1{\null}}
 \def\ShowReservedBoxes{\global\def\FrameSpider##1{##1}}
 \let\HideDisplacementBoxes\HideReservedBoxes  
 \let\ShowDisplacementBoxes\ShowReservedBoxes
 \let\HideFigureFrames\HideReservedBoxes
 \let\ShowFigureFrames\ShowReservedBoxes
  \ShowDisplacementBoxes
 
 \def\hSlide#1{\advance\XSlide@ by #1}
 \def\vSlide#1{\advance\YSlide@ by #1}
 
  \def\SetInkShift@{%
            \advance\XShift@ by -\LT@
            \advance\XShift@ by \XSlide@
            \advance\YShift@ by -\BT@
            \advance\YShift@ by -\YSlide@
             }
  \def\InkShift@#1{\Shifted@{\Scrunched{#1}}}
 
   %
  \def\CleanRegisters@{%
      \globaldefs=1\relax
        \XShift@=\z@\YShift@=\z@\XSlide@=\z@\YSlide@=\z@
        \TT@=\z@\LT@=\z@\BT@=\z@\RT@=\z@
      \globaldefs=0\relax}

 
 \def\SetTexturesEPSFSpecial{\PSOriginfalse
  \gdef\EPSFSpecial##1##2{\relax
    \edef\specialthis{##2}%
    \SPLIT@0.@\specialthis.@\relax
    \special{illustration ##1 scaled
                        \the\Initialtoks@}}}
 
  \def\SetUnixCoopEPSFSpecial{\PSOrigintrue 
   \gdef\EPSFSpecial##1##2{%
      \dimen4=##2pt
      \divide\dimen4 by 1000\relax
      \Real{\dimen4}
      \edef\Aux@{\the\Realtoks}%
      \includegraphics{##1\space}}}

  \def\SetBechtolsheimEPSFSpecial{\PSOrigintrue 
   \gdef\EPSFSpecial##1##2{%
      \dimen4=##2pt
      \divide\dimen4 by 1000\relax
      \Real{\dimen4}
      \edef\Aux@{\the\Realtoks}%
      \special{ps: psfiginit}%
      \special{ps: literal 1 1 0 0 1 1 startTexFig
           \the\mag\space 1000 div \Aux@\space mul 
           \the\mag\space 1000 div \Aux@\space mul scale}%
      \special{ps: include  ##1}%
      \special{ps: literal endTexFig}%
        }}

  \def\SetLisEPSFSpecial{\PSOrigintrue 
   \gdef\EPSFSpecial##1##2{%
      \dimen4=##2pt
      \divide\dimen4 by 1000\relax
      \Real{\dimen4}
      \edef\Aux@{\the\Realtoks}%
      \special{pstext="1 1 0 0 1 1 startTexFig\space
           \the\mag\space 1000 div \Aux@\space mul 
           \the\mag\space 1000 div \Aux@\space mul scale}%
      \includegraphics{##1}%
      \special{pstext=endTexFig}%
        }}

  \def\SetRokickiEPSFSpecial{\PSOrigintrue 
   \gdef\EPSFSpecial##1##2{%
      \dimen4=##2pt
      \divide\dimen4 by 10\relax
      \Real{\dimen4}
      \edef\Aux@{\the\Realtoks}%
      \includegraphics{##1}}}

  \def\SetInlineRokickiEPSFSpecial{\PSOrigintrue 
   \gdef\EPSFSpecial##1##2{%
      \dimen4=##2pt
      \divide\dimen4 by 1000\relax
      \Real{\dimen4}
      \edef\Aux@{\the\Realtoks}%
      \special{ps::[begin] 1 1 0 0 1 1 startTexFig\space
           \the\mag\space 1000 div \Aux@\space mul 
           \the\mag\space 1000 div \Aux@\space mul scale}%
      \special{ps: plotfile ##1}%
      \special{ps::[end] endTexFig}%
        }}

  \def\SetOzTeXEPSFSpecial{\PSOriginfalse 
  \gdef\EPSFSpecial##1##2{
     \special{##1\space 
       ##2 1000 div \the\mag\space 1000 div mul
       ##2 1000 div \the\mag\space 1000 div mul scale
       \the\LLXtoks@\space neg 
       \the\LLYtoks@\space neg translate
             }}} 
  
 \def\SetOzTeXPreviewedEPSFSpecial{\PSOrigintrue
 \gdef\EPSFSpecial##1##2{%
 \dimen4=##2pt
 \divide\dimen4 by 1000\relax
 \Real{\dimen4}
 \edef\Aux@{\the\Realtoks}
 \special{epsf="##1"\space scale=\Aux@}%
 }} 

  \let\SetPSprintEPSFSpecial\SetOzTeXEPSFSpecial
  \let\SetPsprintEPSFSpecial\SetOzTeXEPSFSpecial

 \def\SetArborEPSFSpecial{\PSOriginfalse 
   \gdef\EPSFSpecial##1##2{%
     \edef\specialthis{##2}%
     \SPLIT@0.@\specialthis.@\relax 
     \special{ps: epsfile ##1\space \the\Initialtoks@}}}

 \def\SetClarkEPSFSpecial{\PSOriginfalse 
   \gdef\EPSFSpecial##1##2{%
     \Rescale {\Wd@@}{##2pt}{1000pt}%
     \Rescale {\Ht@@}{##2pt}{1000pt}%
     \special{dvitops: import 
           ##1\space\the\Wd@@\space\the\Ht@@}}}

  \let\SetDVIPSONEEPSFSpecial\SetUnixCoopEPSFSpecial
  \let\SetDVIPSoneEPSFSpecial\SetUnixCoopEPSFSpecial

  \def\SetBeebeEPSFSpecial{
   \PSOriginfalse%
   \gdef\EPSFSpecial##1##2{\relax
    \special{language "PS"
      literal "##2 1000 div ##2 1000 div scale
      position = "bottom left",
      include "##1"}}}
  \let\SetDVIALWEPSFSpecial\SetBeebeEPSFSpecial

  \def\SetNorthlakeEPSFSpecial{\PSOrigintrue
   \gdef\EPSFSpecial##1##2{%
     \edef\specialthis{##2}%
     \SPLIT@0.@\specialthis.@\relax 
     \special{insert ##1,magnification=\the\Initialtoks@}}}

 \def\SetStandardEPSFSpecial{%
   \gdef\EPSFSpecial##1##2{%
     \ms@g{}
     \ms@g{%
       !!! Sorry! There is still no standard for \string%
       \special\ EPSF integration !!!}%
     \ms@g{%
      --- So you will have to identify your driver using a command}%
     \ms@g{%
      --- of the form \string\Set...EPSFSpecial, in order to get}%
     \ms@g{%
      --- your graphics to print.  See BoxedEPS.doc.}%
     \ms@g{}
     \KillEPSFSpecial
     }}

  \def\KillEPSFSpecial{\gdef\EPSFSpecial##1##2{}}

  \SetStandardEPSFSpecial 
 
 \let\wlog\wlog@ld 

 \catcode`\:=\C@tColon
 \catcode`\;=\C@tSemicolon
 \catcode`\?=\C@tQmark
 \catcode`\!=\C@tEmark

 \catcode`\@=\CatAt

 %
 %
 %
 %
 %
\let\epsffile\EPSFfile
 
  \def\FigureTitle#1{\medskip%
        \centerline{\Figurefont #1\unskip}%
        }  
    
  \catcode`\@=11
     
 \let\Mas\relax\let\mas\relax
 \let\Sam\relax\let\sam\relax
 \let\HideEdStuff\relax
 \let\PrintEdStuff\relax
 \let\KillEdStuff\relax
 \let\ShowEdStuff\relax
 \let\change\relax
 \let\beginchange\relax
 \let\endchange\relax
 \let\cbar\relax
 \let\ccbar\relax
 \let\Ninepoint\relax
 \let\ninepoint\relax
 \def\StdMathsurround{\dimen0}
 \let\LoadCMFonts\relax 
 \let\LoadPSFonts\relax 
 \let\LoadNinepoint\relax 
 \let\lBr[ 
 \let\rBr] 
 \let\SetdeGEPSFSpecial\relax
 \def\cal{\fam\tw@}

 \def\SetAuthorHead#1{}
 \def\SetTitleHead#1{}
 \def\ednote#1{}
 \def\Ednote#1{}

  \ifx\amstexloaded@\relax 
          \def\textonlyfont@#1#2{\def#1{#2}}
          \def\textfont@#1#2{\def#1{#2}}
  \else 
   \fi 
   
   \let\wlog\wlog@ld 
   \catcode`\@=\CatAt  

 

 \input amssym.def \input amssym.tex

   \def \C {\Bbb C} \def \R {\Bbb R} 
\def \M{\Bbb M}
\def \Z{\Bbb Z}
\def \C{\Bbb C}
\def \R{\Bbb R}
\def \Q{\Bbb Q}
\def \N{\Bbb N}
\def \l{\lambda}
\def \V{V^{\natural}}
\def \wt{{\rm wt}}
\def \tr{{\rm tr}}
\def \Res{{\rm Res}}
\def \End{{\rm End}}
\def \Aut{{\rm Aut}}
\def \mod{{\rm mod}}
\def \Hom{{\rm Hom}}
\def \<{\langle} 
\def \>{\rangle} 
\def \w{\omega}
\def \o{\omega}
\def \t{\tau }
\def \a{\alpha }
\def \b{\beta}
\def \e{\epsilon }
\def \L{\Lambda}
\def \la{\lambda }
\def \om{\omega }
\def \O{\Omega}
\def \voa{vertex operator algebra\ }
\def \voas{vertex operator algebras\ }
\def \p{\partial}
\def \D{\Delta}
\def \P{\Phi}

   \SetRokickiEPSFSpecial  
   \ShowFigureFrames 

 
 \H           Associative subalgebras of the Griess algebra and related 
topics
 
       C. Dong,\ \ H. Li,\ \ G. Mason,\ \ S. P. Norton

\HH 1. Introduction

In [DMZ] it was shown that the moonshine module $\V$ contains a sub
\voa isomorphic to the tensor product $L({1\over 2},0)^{\otimes 48}$
where $L({1\over 2},0)$ is the \voa associated to highest weight
unitary representation for the Virasoro algebra with central charge
${1\over2}.$ This containment turns out to be fundamental, because it
allows us to deduce properties of $\V$ from those of $L({1\over
2},0),$ which are very much easier to discern
(loc. cit.). This approach has been used to prove, for example,  that $\V$ is
holomorphic [D], and to construct twisted sectors and intertwining
algebras for $\V$ [DLM], [H]. One can even base a simplified approach
to the existence of $\V$ on $L({1\over 2},0)^{\otimes 48}$
(cf. Miyamoto's lecture at this conference).

The reason why $L({1\over 2},0)$ is so attractive is that it is the
first non-trivial {\it discrete series} representation of the Virasoro
algebra and all modules and fusion rules are known for this family of
vertex operator algebras. It is therefore a natural question to ask if
$\V$ contains other similar tensor products of discrete series
representations, and more generally to ask which discrete series
representations can be generated by idempotents in the Griess algebra
$B.$ We do not have complete answers to these questions, but we will
take an approach which allows us, for example, to exhibit in a fairly
painless way a sub \voa of $\V$ of the shape $$\left(L({1\over
2},0)\bigotimes L({7\over 10},0)\bigotimes L({4\over
5},0)\right)^{\otimes12}
\eqno(1.1)$$
as well as another (of central charge less than $24$) which is the tensor
product of the first 24 members of the discrete series. It turns out
that the key idea is to relate these questions to the theory of root
systems and Niemeier lattices.

These issues are naturally related to the question of describing the
maximal associative subalgebras of the Griess algebra. In this form,
these question were first studied in [MN]. We show that each Niemeier
lattice (and its attendant root system) determines (in many ways)
certain maximal associative subalgebras of $B.$ For example, (1.1) is
associated to the Niemeier lattice of type $A_2^{12},$ while
$L({1\over 2},0)^{\otimes48}$ corresponds to $A_1^{24}.$

G.M. would like to thank Gerald Hoehn for stimulating conversations
and Jeff Harvey for pointing out not only the occurrence of the
parafermion in Theorem 2.7, but also that he (together with Lance Dixon) had
done a number of the calculations in Section 2 independently (and much
earlier).  We gratefully acknowledge the following partial support:
C.D.: a faculty research funds granted by the University of
California, Santa Cruz and DMS-9303374; G.M.: NSF grant DMS-9401272.

\HH 2. Algebras and root systems

We consider a simple (i.e., irreducible) root system $\Phi$ of type
$A,D,E.$ Let $\P^+$ be the set of positive roots and let $l\geq1$ be the 
rank of $\P,$ $N=|\P^+|,$ and $h$ the Coxeter number of $\P.$ It is well-known
that we have 
$$2N=lh.\eqno(2.1)$$

Any $\a\in\P^+$ determines a partition of $\P^+,$ namely
$$\eqalign{&\D_0(\a)=\{\a\}\cr
& \D_1(\a)=\{\beta\in\P^+|(\a,\b)\ne 0, \beta\ne\a\}\cr
& \D_2(\a)=\{\beta\in\P^+|(\a,\b)=0, \beta\ne\a\}.\cr}\eqno(2.2)$$
Here, $(\cdot,\cdot)$ denotes the usual inner product associated with 
$\P,$ normalized so that $(\a,\a)=2$ for $\a\in\P.$ We often
write $\a\sim \b$ in case $\b\in\D_1(\a);$ of course $\sim$ is a
symmetric relation. 

We will define a certain $\Q$-algebra $A=A(\P).$ Additively it is a
free abelian group with a distinguished basis consisting of elements
$t(\a),$ $u(\a)$ with $\a\in\P^+.$ Thus $A$ has rank $2N.$ To define
multiplication, note that if $\a\sim\b$ with $\a,\b\in\P^+$ then there
is a unique $\gamma \in\P^+$ such that $\a\sim\gamma\sim\b.$ We use this
observation to define $$t(\a)t(\b)=\left\{\matrix{
t(\a)+t(\b)-t(\gamma), & {\rm if}\ \a\sim\b\sim\gamma\sim\a\cr 0, & {\rm
if}\ \b\in\D_2(\a)\cr}\right.$$ 
$$u(\a)u(\b)=\left\{\matrix{
u(\a)+u(\b)-t(\gamma), & {\rm if}\ \a\sim\b\sim\gamma\sim\a\cr 0, & {\rm
if}\ \b\in\D_2(\a)\cr}\right.$$
$$u(\a)t(\b)=t(\b)u(\a)=\left\{\matrix{ u(\a)+t(\b)-u(\gamma), & {\rm
if}\ \a\sim\b\sim\gamma\sim\a\cr 0, & {\rm if}\
\b\in\D_2(\a)\cr}\right.$$ $$ t(\a)^2=8t(\a), u(\a)^2=8u(\a).$$

Obviously, $A=A(\P)$ is a commutative (but in general non-associative)
algebra.  The motivation for introducing this algebra will be
explained later.

The following is well-known (cf. Chap. IV, 1.11, Proposition 32 of [B]).

\th Lemma 2.1 

If $\a\in\P^+$ then $|\D_1(\a)|=2h-4.$ 
\endth

We use this to prove

\th  Proposition 2.2

$A(\P)$ has an identity element, namely 
$$\delta={1\over 4h}
\sum_{\a\in\P^+}(t(\a)+u(\a)).$$
\endth

\pf Proof

If $\b\in\P^+$ then $t(\b)\delta$ is equal to
$$\eqalign{
&\ \   \ \ {1\over 4h}\left(t(\b)^2+t(\b)\sum_{\a\in\D_1(\b)}t(\a)
+t(\b)\sum_{\a\in\D_1(\b)}u(\a)\right)\cr
& ={1\over 4h}\left(8t(\b)+\sum_{\a\in\D_1(\b)}(t(\b)+t(\a)-t(\gamma(\a,\b)))
\right.\cr
&\ \ \ \ \ \ \ \left.+\sum_{\a\in\D_1(\b)}(t(\b)+u(\a)-u(\gamma(\a,\b)))\right)\cr}\eqno(2.3)$$
where $\gamma(\a,\b)$  is the element of $\P^+$  determined by $\a$ and
$\b$ whenever $\a\sim \b.$ But $\gamma(\a,\b)$ ranges over $\D_1(\b)$ as 
$\a$ does, so all  terms in (2.3) cancel except for the
$t(\b)$'s. Lemma 2.1 tells us that (2.3) is thus equal
to ${1\over 4h} (8t(\b)+(2h-4)t(\b)+(2h-4)t(\b))=t(\b),$ that is 
$t(\b)\delta=t(\b).$ Similarly $u(\b)\delta=u(\b),$ and the proposition is proved. \qed

We observe that the $\Q$-span of the $t(\a)$ for
$\a\in\P^+$ is a subalgebra of $A,$ which we denote by $T(\P).$ The same proof
shows 

\th Lemma 2.3 

$T(\P)$ has an identity, namely
 $$\epsilon={1\over 4h}
\sum_{\a\in\P^+}t(\a).$$ 
\endth

Now introduce a symmetric bilinear form $\<,\>$ on $A(\P)$ as follows:
$$\<t(\a),t(\a)\>=\<u(\a),u(\a)\>=4$$
$$ \<t(\a),t(\b)\>=\<u(\a),u(\b)\>=\left\{\matrix{
1/2, & \a\sim\b\cr
0, & \a\in\D_2(\b)\cr}\right.$$
$$ \<t(\a),u(\b)\>=\left\{\matrix{
1/2, & \a\sim\b\cr
0, & \a\in\D_0(\b)\cup\D_2(\b).\cr}\right.$$
We remark that $\<,\>$ is not necessarily non-degenerate (see below).

For an element $a\in A$ we define 
$$c(a)=8\<a,a\>.\eqno(2.4)$$
We may refer to $c(a)$ the {\it central charge} of $a,$
though usually we reserve this term for the case that $a$ is an idempotent
of $A$ such as $\delta$ or $\epsilon.$ In these cases we prove  

\th Lemma 2.4 

We have 
$$ c(\delta)=l,\ \ c(\epsilon)={lh\over h+2}.$$ 
\endth

\pf Proof

We have 
$$\eqalign{
&c(\delta)={8\over (4h)^2}\sum_{\a\in\P^+}(\<t(\a),t(\a)\> +
\<u(\a),u(\a)\>)\cr
&\ \ \ \ \ \  \ \ +{8\over (4h)^2}
\sum_{\scriptstyle\a\in\P^+\atop\scriptstyle\b\in\D_1(\a)}(\<t(\a),t(\b)\> +
\<u(\a),u(\b)\>+\<t(\a),u(\b)\> +
\<u(\a),t(\b)\>)\cr 
&\ \ \ \ \ ={1\over 2h^2}(8N+N(2h-4)4)\cr
&\ \ \ \ \ ={2N\over h}.\cr}$$

Now use (2.1) to see that $c(\delta)=l.$ The calculation of $c(\epsilon)$ is 
similar. \qed

We now consider the possibility of decomposing the identity $\delta$ of $A$ 
(or the identity $\epsilon$ of $T$) into a sum of idempotents
$\delta=e_1+e_2+\cdots +e_k$ such that $\<e_i,e_j\>=0$ if
$i\ne j$ and $e_ie_j=\delta_{i,j}e_i.$ Such a decomposition corresponds to a 
particular kind of associative subalgebra of $A$ isomorphic to
a direct sum of $k$ copies of $\Q.$ One way to find such a decomposition
is as follows: locate within the root system $\P$ a subsystem $\P',$
so that $A(\P')$ is a subalgebra of $A(\P).$ Then decompose $\delta$
as $\delta'+(\delta-\delta')$ where $\delta'$ is the identity of $A(\P').$ This
has the desired properties (with $k=2$). For large values of $k$ one 
can iterate this procedure, considering chains of root systems $\P
\supset \P'\supset\P''\cdots.$

We illustrate this procedure in the case that $\P$ is of type
$A_l$ -- the case of the most interest to us. Having fixed a system
$\pi=\{\a_1,...,\a_l\}$ of simple roots in $\P,$ we let $\P_i$ be the
sub-system whose simple roots consist of $\a_1,...,\a_i,$ with $A_i$ and $T_i$
the corresponding algebras. Thus $A=A_l$ and $T=T_l.$

We begin by writing $\delta=\epsilon+(\delta-\epsilon)$ and verifying
that indeed $\<\epsilon,\delta-\epsilon\>=0,$ which we leave to the reader. Then we find that 
$$c( \delta-\epsilon)=c(\delta)-(\epsilon)={2l\over h+2}.$$
As $\P$ is of type $A_l$ we have $h=l+1,$ so that
$$c(\delta-\epsilon)={2l\over l+3}.\eqno(2.5)$$

Now consider the containment $T_{l-1}\subset T_l$ and let $\epsilon'$ be
the identity of $T_{l-1},$ with $l',N'$ having the obvious meaning.
\th Lemma 2.5 

We have
$\<\epsilon-\epsilon',\epsilon'\>=0.$
\endth

 \pf Proof

If $h'=h-1$ is the Coxeter number of $\P_{l-1}$ then by Lemma 2.4 we have
$\<\e,\e'\>={(l-1)h'\over 8(h'+2)}.$ On the other hand we have
$$ \<\e,\e'\>={1\over (2h+4)(2h'+4)}\sum_{\a\in\P^+_l}\sum_{\b\in\P^+_{l-1}}
\<t(\a),t(\b)\>.\eqno(2.6)$$
If $\a\in \P^+_{l-1}$ then $\sum_{\b\in\P^+_{l-1}}
\<t(\a),t(\b)\>=4+{2h'-4\over 2}=h'+2.$ If
$\a\in\P^+_l\backslash \P^+_{l-1}$ then $\sum_{\b\in\P^+_{l-1}}
\<t(\a),t(\b)\>={l-1\over 2}$ since $|\D_1(\a)\cap \P^+_{l-1}|=l-1$ for
each $\a.$ Since $|\P^+_l\backslash \P^+_{l-1}|=l,$ we find that
(2.6) yields 
$$\eqalign{\<\e,\e'\>&={1\over (2h+4)(2h'+4)}\left(N'(h'+2)+{l(l-1)\over 2}\right)\cr
&={1\over 4(h'+2)(h'+3)}\left({h'l'(h'+2)\over 2}+{l(l-1)\over 2}\right)\cr
&={h'(h'+3)(l-1)\over 8(h'+2)(h'+3)}\cr
&={h'(l-1)\over 8(h'+2)}\cr}.$$
So indeed $\<\e,\e'\>=\<e',e'\>,$ as desired. \qed

\th Lemma 2.6

We have $c(e-e')=1-\displaystyle{{6\over (l+2)(l+3)}}.$
\endth

\pf Proof

After Lemma 2.5 we have $c(\e-\e')=c(\e)-c(\e').$ Now use Lemma 2.4 to see that

$$\eqalign{c(\e-\e')&={lh\over h+2}- {l'h'\over h'+2}\cr
&={l(l+1)h\over l+3}- {l(l-1)\over l+2}\cr
&=1-{6\over (l+2)(l+3)}.\cr}$$

If we iterate this procedure, we obtain the following result.

\th Theorem 2.7

Suppose that $A=A(\P)$ corresponding to the root system $\P$ of type $A_l.$ Then the identity $\delta$ of $A$ can be decomposed into a sum of $l+1$
idempotents $e_1,e_2,...,e_{l+1}$ satisfying the following

(i) $\delta=e_1+e_2+\cdots e_{l+1}$

(ii) $e_ie_j=\delta_{ij}e_i$

(iii) $\<e_i,e_j\>=0$ if $i\ne j.$

(iv) $c(e_i)=\displaystyle{1-{6\over (i+2)(i+3)}},$ $1\leq i\leq l.$

(v) $c(e_{l+1})=\displaystyle{{2l\over l+3}.}$
\endth

\rk Remark

The significance of (iv), of course, is that the series of values
$1-{6\over (i+2)(i+3)},$ $i=1,2,...$ correspond to the central 
charges of the discrete series representations of the Virasoro algebra.
These are the unitary highest weight representations of the 
Virasoro algebra with central charge $c$ satisfying $0<c<1.$ 
These values were first identified in [FQS], and the unitarity was proved 
in [GKO].
\endrk

We owe to Jeff Harvey the observation that the central charge
${2l\over l+3}$ also has some significance, namely it corresponds to
the central charge of the parafermion algebras [ZF], [DL]. Both the
discrete series and parafermionic representations of the Virasoro
algebra arise from the GKO ``coset construction,'' indeed our
arguments leading to Theorem 2.7 amount to a coset construction at the
level of root systems. We pursue these ideas in the following 
sections. 

It goes without saying that (i)-(iii) and a modification of (iv) and
(v) also hold if $A_l$ is replaced by root systems of type
$D$ or $E.$

\HH 3. The \voa $V^+_{\sqrt{2}R}$

Suppose that $L$ is a positive-definite even lattice. The \voa $V_L$
associated with $L$ is one of the most fundamental examples of a
vertex operator algebra. We refer the reader to [FLM] for the
construction of $V_L.$ We wish to focus on the case in which $R$ is
the root lattice corresponding to a simple root system $\P$ as in
Section 2, and where $L=\sqrt{2}R$ (later, we will relax our
conditions to allow $R$ to be semi-simple). Thus $L$ is indeed a
positive-definite even lattice which has, in addition, no vectors of
squared length 2.

The lattice $L$ has a canonical automorphism $t$ of order 2, namely
the reflection automorphism $t:x\mapsto -x$ for $x\in L.$ Then $t$
lifts to a canonical automorphism of $V_L$ (loc. cit.), and we denote
by $V^+_L$ the sub \voa of $t$-fixed points of $V_L.$ In the case that
$L=\sqrt{2}R$, $V_L^+$ has no vectors of weight 1 owing to the absence
of vectors in $L$ of squared length 2. Thus $$V^+_{\sqrt{2}R}=\C{\bf
1}\oplus B^+\oplus\cdots \eqno(3.1)$$ where in (3.1) we have set
$B^+=(V^+_{\sqrt{2}R})_2$, i.e., $B^+$ is the space of vectors in
$V^+_{\sqrt{2}R}$ of weight 2.

There is a canonical structure of commutative, non-associative algebra
on $B^+,$ namely
$$ab=a_1b, \ a,b\in B^+$$
where the vertex operator for $a$ is given by
$$Y(a,z)=\sum_{n\in\Z}a_nz^{-n-1}.$$
Similarly, there is a non-degenerate, invariant inner product $\<\cdot,\cdot\>$
on $B^+$ given by
$$\<a,b\>{\bf 1}=a_3b.$$

Conveniently, the algebra structure of $B^+$ has been written down
in Theorem 8.9.5 and Remark 8.9.7 of [FLM]. To describe it, let 
$H=\Q\bigotimes\sqrt{2}R.$ Then 
$$ B^+=S^2(H)\oplus \oplus_{\a\in\P^+}\Q x_{\a}$$
where 
$$ x_{\a}=e^{\sqrt{2}\a}+e^{-\sqrt{2}\a}$$
and $e^{\sqrt{2}\a}$ is the standard notation for the element
in the group algebra $\C[L]$ corresponding to $\sqrt{2}\a.$
(Note that in this case $\hat L$ is a direct product of $L$ and $\<\pm 1\>.$)
The relations are as follows:
$$\eqalign{h^2\cdot k^2&=4(h,k)hk,\ \ h,k\in H\cr
h^2\cdot x_{\a}&=2(h,\a)^2x_{\a}\cr
x_{\a}\cdot x_{\b}&=\left\{\matrix{0, & \b\in\D_2(\a)\cr
x_{\gamma},& \a\sim\b\sim\gamma\sim \a\cr
2\a^2, & \a=\b.\cr}\right.\cr}$$
Moreover
$$\eqalign{\<h^2,k^2\>&=2(h,k)^2\cr
\<h^2,x_{\a}\>&=0\cr
\<x_{\a},x_{\b}\>&=\left\{\matrix{0,&\a\ne\b\cr
2,& \a=\b.\cr}\right.\cr}$$

\th Theorem 3.1

Let $\P$ be a simple root system of type $A, D,E,$ let $A=A(\P)$ be as in 
Section 2, and let $B^+$  be as above. There there is an isometric surjection
of algebras $A\to B^+$ given by
$$\eqalign{\phi: &t(\a)\mapsto {1\over 2}\a^2-x_{\a}\cr
& u(\a)\mapsto {1\over 2}\a^2+x_{\a}.\cr}$$
\endth

\th Corollary 3.2

If $\P$ is of type $A_l$ then the map $\P:A\to B^+$ is an isometric isomorphism
of algebras
\endth

Note that $\dim B^+={l(l+1)\over 2}+N.$ If $\P$ is
of type $A_l$ then $N={l(l+1)\over 2},$ so that $\dim B^+=2N=\dim A.$ 
So in this case, Theorem 3.1 implies the corollary.

 \rk Remark 

In the other cases, i.e., $\P$ of type $D$ or $E,$ we similarly see that
the surjection $A\to B^+$ is {\it not} an isomorphism. Since it is an
isometry (by the theorem), the kernel is 
precisely the radical of the form $\<,\>$ on $A.$
\endrk

\pf Proof of Theorem

This is quite straightforward. For example if $\a\sim\b\sim\gamma\sim\a$ then
$$\phi(u(\a))\phi(t(\b))=({1\over 2}\a^2+x_{\a})({1\over 2}\b^2-x_{\b})
=(\a,\b)\a\b-x_{\gamma}+{1\over 2}(\b^2 x_{\a}-\a^2 x_{\b}).\eqno(3.2)$$
In this case we have $(\a,\b)=-1$ and $\a\b={1\over 2}(\gamma^2-\a^2-\b^2)$ 
(since $\a+\b+\gamma=0$). So (3.2) reads 
$${1\over 2}(\a^2+\b^2-\gamma^2)-x_{\gamma}+x_{\a}-x_{\b}
=\phi(u(\a))+\phi(t(\b))-\phi(u(\gamma))$$
as required. The other relations which show that
$\phi$ is an algebra morphism are proved similarly. It is
clear that $\phi$ is an epimorphism, so we only need check that it is also
an isometry. This is also easy;
for example
$$\<\phi(u(\a)),\phi(t(\b))\>=\<{1\over 2}\a^2+x_{\a},{1\over 2}\b^2-x_{\b}\>
={1\over 2}(\a,\b)^2-2\delta_{\a,\b}=\<t(\a),u(\b)\>.$$
The other cases are similar. \qed

Denote by $L(c,0)$ the vertex operator algebra associated to the 
Virasoro algebra with central charge $c.$ We have

\th Theorem 3.3

Let $\P$ be the simple root system of type $A_l,$ with $R$ the
corresponding root lattice. Then the \voa $V^+_{\sqrt{2}R}$ contains a
sub \voa isomorphic to  a tensor product 
$$\bigotimes_{i=1}^{l+1}L(c_i,0)$$
where $c_i=1-{6\over (i+2)(i+3)}$ for $1\leq i\leq l$ and 
$c_{l+1}={2l\over l+3}.$

\pf Proof

By a simple calculation, Miyamoto has shown (Theorem 4.1 of [M]) that if
$e$ is an idempotent in $B^+$ then $2e$ is such that the components of
the vertex operator $Y(2e,z)$ generate a copy of the Virasoro algebra
with central charge $c(e)$ in the sense of (2.4).
Moreover, orthogonal idempotents $e_i,e_j$ satisfying $e_ie_j=0$
are such that the corresponding Virasoro algebras commute, that is they
generate a tensor product of the two algebras (cf. Lemma 5.1 of (loc. cit.)
and Lemma 2.4 of [DLM]). Now the theorem follows from Theorem 2.7.
\qed

\rk Remark

The theorem has an obvious extension to semi-simple root systems, all of 
whose simple components are of type $A_l$ (varying $l$). 
\endrk

\HH 4. Niemeier lattices and the Moonshine module

Recall that a {\it Niemeier lattice}\footnote*{\rm Somewhat 
unconventionally,
according to this terminology the Leech lattice falls under the rubric of 
Niemeier lattice.} is a self-dual, even, positive-definite
lattice $L$ of rank 24. The set of vectors $L_2$ of squared-length 2 forms
a root system which is either empty (in which case $L$ is the Leech
lattice) or else itself has rank 24. In this latter case the root system $L_2$ is such that 
its simple components have a common Coxeter number. We call such a root system
a Niemeier root system; then the map
$$L\to L_2$$
sets up a bijection between Niemeier lattices and Niemeier root systems 
(possibly empty). See [Nie] and [V] for background.

We now discuss the following result, which was hinted at in [Nor, P305]
and which extends the results of [Cur]. 

\th Theorem 4.1

Let $L$ be any Niemeier lattice. Then there is at least one (and in
general several) isometric embedding $\sqrt{2}L\to \L,$ where $\L$ is
the Leech lattice.
\endth

We want to know how many sublattices of the Leech lattice $\Lambda$ there are 
which are versions of each Niemeier lattice, including the Leech lattice itself,
rescaled by a factor of $2$. For each type of lattice the mass of sublattices of
that type is shown in Table 1; the number of sublattices can be obtained from 
this by multiplying by the order of the Conway group $Co_1$, which is 
$2^{21}3^95^47^211.13.23$.

We now explain how the results in Table 1 may be obtained. First, it
is easily seen (see Lemma 4.3 below) that if $N=\sqrt{2}L\subset \L$
with $L$ a Niemeier lattice then $N$ corresponds to a Lagrangian
(i.e., maximal totally isotropic) subspace of $\L/2\L$, under the
orthogonal form whose value on a vector is 0 or 1 according 
as the squared length of the corresponding Leech lattice vector is or is not divisible
by 4. This form is of type $+,$ so that a Lagrangian subspace has dimension 12.

Now the number of extensions of any isotropic subspace to a Lagrangian one
is known. To be precise, if the space has codimension $n$ (i.e.\ dimension 
$12-n$) the number will be
$$\prod^{n-1}_{i=0}(2^i+1)\eqno(4.1)$$
this being interpreted as $1$ if $n=0$. The problem now is to determine how to 
use this information to work out the number of Lagrangian subspaces of each 
type.

We do this by induction on the codimension. Let us call an isotropic subspace 
{\sl connected} if it is generated by its minimal (norm $4$) vectors, and these 
cannot be split into components that generate disjoint subspaces. It is well 
known that connected spaces correspond to Dynkin diagrams of type $A_n$ 
($n\ge 1$), $D_n$ ($n\ge 4$) or $E_n$ ($n=6,7,8$). We classify all subspaces of 
each of these types and, for each of them, determine how many extensions it has 
to connected spaces of the next higher dimension (or, equivalently, how many 
subspaces of the next lower dimension belong to each type.

This gives us the information we need. For by induction we already know how many
Lagrangian spaces contain any connected subspace of higher dimension, and
we know by (4.1) the total number of extensions. The difference between these 
will, therefore, be the number of extensions which contain a subspace of the 
relevant type as a connected component.

To take some examples, any Lagrangian subspace in which $A_8$ appears as 
a connected component must correspond to a Niemeier lattice of type $A_8^3$. If 
$D_6$ appears as a component then it will be of type $D_6^4$ or $A_9^2D_6$; we 
can distinguish between these by counting the number of extensions of an $A_9$.

\bigbreak
\centerline{\bf Table 1}
\bigbreak
$$\vbox{
\offinterlineskip
\halign{\strut\vrule#&\hfil\quad$#$\quad&\vrule#&\quad$#$\quad\hfil&#\vrule\cr
\noalign{\hrule}
&\Lambda&&153715/123771648&\cr
&A_1^{24}&&141985575/58032128&\cr
&A_2^{12}&&469525/19008&\cr
&A_3^8&&3077275/86016&\cr
&A_4^6&&39821/2560&\cr
&A_5^4D_4&&24653/5120&\cr
&D_4^6&&58025/524288&\cr
&A_6^4&&16813/20160&\cr
&A_7^2D_5^2&&6087/20480&\cr
&A_8^3&&1765/64512&\cr
&A_9^2D_6&&4037/276480&\cr
&D_6^4&&791/294912&\cr
&A_{11}D_7E_6&&67/46080&\cr
&E_6^4&&13/184320&\cr
&A_{12}^2&&61/366080&\cr
&D_8^3&&575/14680064&\cr 
&A_{15}D_9&&41/3870720&\cr
&A_{17}E_7&&1/645120&\cr
&D_{10}E_7^2&&11/5160960&\cr
&D_{12}^2&&19/259522560&\cr
&A_{24}&&1/244823040&\cr
&D_{16}E_8&&1/660602880&\cr
&E_8^3&&1/1981808640&\cr
&D_{24}&&1/501397585920&\cr
\noalign{\hrule}}}$$
\bigbreak

The computations were done in {\bf MAPLE.}

Table 2 shows the complete list of connected isotropic subspaces. The first 
column shows the corresponding symbol, with Greek letters used to distinguish 
subspaces with the same Dynkin diagram. The second column shows the dimension. 
The third column gives the maximal intersection with the lattice of type 
$A_{23}$ generated by vectors of type $(4,-4,0^{22})$ (parametrized by the 
corresponding subset of the coordinate positions, on which the Conway group 
determines a Gloay code), plus the number of ``spare'' vectors. The symbols 
$s_{12}^+$, $s_n$ and $u_n$ are taken from [C1]: the first denotes the union of 
two special octads which intersect in a tetrad; the second a set of cardinality 
$n$ in an ascending chain containing a special octad and the complement of 
another such, and {\sl not} containing an $s_{12}^+$; and the third a set of the
right size in an ascending chain linking a non-special hexad (i.e., {\sl not} an
$s_6$), a special dodecad, and the complement of a non-special hexad.

The fourth column shows the structure of the stabilizer of the subspace in 
$Co_1$. The subgroup to the left of the aligned column of dots is the
pointwise stabilizer. (Square brackets denote an unspecified group of the 
relevant order.) The last column gives a list of lattices one dimension lower, 
with the number of extensions and containments of the relevant type. 
For example, the entry ``$30\!:\!2(A_8^{\alpha})$'' in the row beginning 
``$A_9^{\alpha}$'' means that an $A_9^{\alpha}$ contains $2$ $A_8^{\alpha}$'s, 
while an $A_8^{\alpha}$ extends to $30$ $A_9^{\alpha}$'s.

\centerline{\bf Table 2}
\bigbreak
\halign{$#$\hfil&\quad$#$\hfil&\quad$#$\hfil&\hfil$\;#.$&$#\;$\hfil&$#$\hfil\cr
0&0&&Co_1&1&\cr
A_1&1&s_2&Co_2&1&98280\!:\!1(0)\cr
A_2&2&s_3&U_6(2)&S_3&2300\!:\!3(A_1)\cr
A_3&3&s_4&2^9.M_{21}&S_4&891\!:\!4(A_2)\cr
D_4&4&s_4+1&2^{4+8}.A_5&2^{1+4}.(S_3\times 3)&21\!:\!12(A_3)\cr
A_4&4&s_5&2^{1+8}.A_5&S_5&336\!:\!5(A_3)\cr
D_5&5&s_5+1&2^7.2^4.3&2^4.S_5&60\!:\!5(D_4),5\!:\!16(A_4)\cr
A_5^{\alpha}&5&s_6&2^{1+8}.S_3&S_6&10\!:\!6(A_4)\cr
A_5^{\beta}&5&u_6&2^6.3&S_6&160\!:\!6(A_4)\cr
E_6&6&s_6+1&2^{1+8}.3&U_4(2).2&4\!:\!27(D_5),1\!:\!36(A_5^{\alpha})\cr
D_6^{\alpha}&6&s_6+1&2^6.2^4&2^5.S_6&3\!:\!6(D_5),3\!:\!32(A_5^{\alpha})\cr
D_6^{\beta}&6&u_6+1&2^6.3&2^5.S_6&16\!:\!6(D_5),1\!:\!32(A_5^{\beta})\cr
A_7^{\alpha}&6&s_8&2^{1+8}&S_8&3\!:\!28(A_5^{\alpha})\cr
A_6^{\alpha}&6&u_7&2^5&S_6&96\!:\!1(A_5^{\alpha}),36\!:\!6(A_5^{\beta})\cr
A_6^{\beta}&6&u_6+1&3&S_7&64\!:\!7(A_5^{\beta})\cr
E_7&7&s_8+1&2^{1+8}&S_6(2)&
 3\!:\!28(E_6),2\!:\!63(D_6^{\alpha}),1\!:\!36(A_7^{\alpha})\cr
D_8^{\alpha}&7&s_8+1&2^{1+8}&2^6.A_8&
 1\!:\!28(D_6^{\alpha}),2\!:\!64(A_7^{\alpha})\cr
D_7&7&u_7+1&2^5&2^6.S_6&
 16\!:\!1(D_6^{\alpha}),18\!:\!6(D_6^{\beta}),1\!:\!64(A_6^{\alpha})\cr
A_8^{\alpha}&7&s_9&2^4&A_8&64\!:\!1(A_7^{\alpha}),2\!:\!28(A_6^{\alpha})\cr
A_7^{\beta}&7&u_8&2^4&2^4.S_4&30\!:\!8(A_6^{\alpha})\cr
A_7^{\gamma}&7&u_7+1&1&S_6\times 2&
 32\!:\!2(A_6^{\alpha}),63\!:\!6(A_6^{\beta})\cr
E_8&8&s_8+2&2^{1+8}&O_8^+(2)&1\!:\!120(E_7),1\!:\!135(D_8^{\alpha})\cr
D_9^{\alpha}&8&s_9+1&2^4&2^7.A_8&
 16\!:\!1(D_8^{\alpha}),1\!:\!28(D_7),1\!:\!128(A_8^{\alpha})\cr
D_8^{\beta}&8&u_8+1&2^4&2^7.(2^4S_4)&15\!:\!8(D_7),1\!:\!128(A_7^{\beta})\cr
A_9^{\alpha}&8&s_{10}&2^3&2^4.L_3(2)&
 30\!:\!2(A_8^{\alpha}),8\!:\!28(A_7^{\beta})\cr
A_9^{\gamma}&8&s_9+1&1&S_8&16\!:\!2(A_8^{\alpha}),1\!:\!28(A_7^{\gamma})\cr
A_8^{\beta}&8&u_9&2^3&3^2.2S_4&16\!:\!9(A_7^{\beta})\cr
A_8^{\gamma}&8&u_8+1&1&2^4.S_4&16\!:\!1(A_7^{\beta}),30\!:\!8(A_7^{\gamma})\cr
A_8^{\delta}&8&u_7+2&1&S_3\wr S_3&10\!:\!9(A_7^{\gamma})\cr
D_{10}^{\alpha}&9&s_{10}+1&2^3&2^8.(2^4.L_3(2))&
 15\!:\!2(D_9^{\alpha}),4\!:\!28(D_8^{\beta}),1\!:\!256(A_9^{\alpha})\cr
D_9^{\beta}&9&u_9+1&2^3&2^8.(3^2.2S_4)&
 8\!:\!9(D_8^{\beta}),1\!:\!256(A_8^{\beta})\cr
A_{11}^{\alpha}&9&s_{12}^+&2^2&2^6.3^{1+2}.2^2&
 14\!:\!18(A_9^{\alpha}),8\!:\!64(A_8^{\beta})\cr
A_{10}^{\alpha}&9&s_{10}+1&1&2^4.L_3(2)&
 8\!:\!1(A_9^{\alpha}),30\!:\!2(A_9^{\gamma}),4\!:\!28(A_8^{\gamma})\cr
A_9^{\beta}&9&u_{10}&2^2&S_6.2&6\!:\!10(A_8^{\beta})\cr
A_9^{\delta}&9&u_9+1&1&3^2.2S_4&8\!:\!1(A_8^{\beta}),8\!:\!9(A_8^{\gamma})\cr
A_9^{\epsilon}&9&u_8+2&1&2^5.S_5&1\!:\!10(A_8^{\gamma})\cr
A_9^{\zeta}&9&u_8+2&1&[2^63]&12\!:\!6(A_8^{\gamma}),27\!:\!4(A_8^{\delta})\cr
D_{12}^{\alpha}&10&s_{12}^++1&2^2&2^9.(2^6.3^{1+2}.2^2)&
 7\!:\!18(D_{10}^{\alpha}),4\!:\!64(D_9^{\beta}),1\!:\!512(A_{11}^{\alpha})\cr
D_{10}^{\beta}&10&u_{10}+1&2^2&
 2^9.S_6.2&3\!:\!10(D_9^{\beta}),1\!:\!512(A_9^{\beta})\cr
A_{15}^{\alpha}&10&s_{16}&2&2^4.A_8&
 6\!:\!140(A_{11}^{\alpha}),4\!:\!448(A_9^{\beta})\cr
A_{12}^{\alpha}&10&s_{12}^++1&1&2^6.3^{1+2}.2^2&
 4\!:\!1(A_{11}^{\alpha}),7\!:\!18(A_{10}^{\alpha}),4\!:\!64(A_9^{\delta})\cr
A_{11}^{\beta}&10&u_{12}&2&M_{12}&2\!:\!66(A_9^{\beta})\cr
A_{11}^{\gamma}&10&s_{10}+2&1&[2^93]&
 7\!:\!4(A_{10}^{\alpha}),10\!:\!4(A_9^{\epsilon}),3\!:\!24(A_9^{\zeta})\cr
A_{10}^{\beta}&10&u_{10}+1&1&S_6.2&
 4\!:\!1(A_9^{\beta}),3\!:\!10(A_9^{\delta})\cr
A_{10}^{\gamma}&10&u_9+2&1&[2^43^2]&
 6\!:\!2(A_9^{\delta}),12\!:\!9(A_9^{\zeta})\cr
D_{16}&11&s_{16}+1&2&2^{10}.(2^4.A_8)&3\!:\!140(D_{12}^{\alpha}),
 2\!:\!448(D_{10}^{\beta}),1\!:\!1024(A_{15}^{\alpha})\cr
D_{12}^{\beta}&11&u_{12}+1&2&2^{10}.M_{12}&
 1\!:\!66(D_{10}^{\beta}),1\!:\!1024(A_{11}^{\beta})\cr
A_{23}&11&s_{24}&1&M_{24}&
 2\!:\!759(A_{15}^{\alpha}),2\!:\!2576(A_{11}^{\beta})\cr
A_{16}&11&s_{16}+1&1&2^4.A_8&2\!:\!1(A_{15}^{\alpha}),
 3\!:\!140(A_{12}^{\alpha}),2\!:\!448(A_{10}^{\beta})\cr
A_{13}&11&s_{12}^++2&1&[2^93^2]&
 3\!:\!2(A_{12}^{\alpha}),6\!:\!18(A_{11}^{\gamma}),2\!:\!64(A_{10}^{\gamma})\cr
A_{12}^{\beta}&11&u_{12}+1&1&M_{12}&
 2\!:\!1(A_{11}^{\beta}),1\!:\!66(A_{10}^{\beta})\cr
A_{11}^{\delta}&11&u_{10}+2&1&S_6\times 2&
 2\!:\!2(A_{10}^{\beta}),1\!:\!10(A_{10}^{\gamma})\cr
A_{11}^{\epsilon}&11&u_9+3&1&3^2.2S_4&4\!:\!12(A_{10}^{\gamma})\cr
D_{24}&12&s_{24}+1&1&2^{11}.M_{24}&
 1\!:\!759(D_{16}),1\!:\!2576(D_{12}^{\beta})\cr
A_{24}&12&s_{24}+1&1&M_{24}&1\!:\!1(A_{23}),1\!:\!759(A_{16})\cr
A_{17}&12&s_{16}+2&1&2^4.A_8\times 2&
 1\!:\!2(A_{16}),1\!:\!140(A_{13}),1\!:\!448(A_{11}^{\delta})\cr
A_{15}^{\beta}&12&s_{12}^++4&1&[2^{12}3^3]&
 1\!:\!24(A_{13}),1\!:\!256(A_{11}^{\epsilon})\cr
A_{12}^{\gamma}&12&u_9+4&1&L_3(3)&1\!:\!13(A_{11}^{\epsilon})\cr}
\vfill\eject

Let us fix a Niemeier lattice $L$ {\it not} equal to the Leech lattice,
and let $R=L_2.$ So there are isometric embeddings 
$$\sqrt{2}R\to \sqrt{2}L\to \L$$
and correspondingly there are \voa embeddings 
$$V^+_{\sqrt{2}R}\to V^+_{\sqrt{2}L}\to V^+_{\L}\to V^{\natural}\eqno(4.2)$$
where $V^{\natural}$ is the moonshine module [FLM].

Suppose that $L_2$ is one of the following type:
$A_1^{24}, A_2^{12},$ or $A_{24}.$ Then we may use (4.2) together with
Theorem 3.3 to see that $V^{\natural}$ contains sub \voas of the 
following kind:
 $$L({1\over 2},0)^{\otimes 48}\eqno(4.3)$$
$$\bigotimes_{i=1}^{24}L(c_i,0)\bigotimes L({16\over 9},0)\eqno(4.4)$$
as well as the type in equation (1.1). 

Type (4.3) was already constructed in [DMZ]. Of course,
we get many analogous tensor products by using 
other Niemeier lattices. Type (1.1) is of interest because,
along with (4.3), it seems to be the only tensor product of Virasoro
algebras that one can obtain in this way which is both of central charge
24 and has only discrete series as factors. The factor 
$L({116\over 117},0)$ in (4.4) (corresponding to $i=24$) corresponds to 
the discrete series with largest value of $c$  which we know occurs via 
an idempotent of the Griess algebra. It would be interesting  to know if it 
is indeed the maximal such value of $c.$ 

Of course, the embeddings (4.2) imply embeddings of the corresponding
homogeneous spaces. In particular, if we let $B^+(L_2)$ be the weight
2 subspace of $V^+_{\sqrt{2}R}$ then we get $B^+(L_2)\subset B.$
Then application of Theorem 3.3 (and the remark following it) 
yields embeddings of associative subalgebras into the Griess algebra
$B.$ For example if $L_2$ is of type $A_1^{24},$ corresponding to
type (4.3), we get a maximal associative subalgebra of $B$ of dimension
48. This was first constructed in [MN]. Similarly type (1.1) 
gives a maximal associative subalgebra of $B$ of dimension 36. In the general
case we have

\th Lemma 4.2

Let $L$ be a Niemeier lattice with root system $L_2$ which has $k$ simple
components. Then $L$ determines (in several ways) an associative subalgebra
of $B$ of dimension $24+k.$
\endth

\pf Proof

The procedure of Theorem 2.7 yields an $(l+1)$-dimensional associative
algebra for any simple root system of rank $l,$ and it maps 
isomorphically into the corresponding $B^+$ of theorem 3.1 because 
$\phi$ is an isometry. Thus our Niemeier lattice $L$ affords an associative
algebra  of dimension $\sum_{i=1}^k(1+l_i)$ where $\{l_i\}$
are the ranks  of the simple components of $L_2.$ Since $\sum_{i}l_i=24,$
the lemma follows. \qed

There is another way to look at these matters which involves the
Monster more directly and which was briefly touched on above and in
[Nor]. Namely, let $\bar \L=\L/2\L$ be the Leech lattice mod 2
equipped with the non-degenerate form which was described earlier.
 Then $\bar \L$ has type $+,$ that is the Lagrangian subspaces have rank 12. If
$\bar U\subset \bar \L$ is such a space then the full inverse image $U\subset \L$ has index $2^{12},$ so that ${1\over\sqrt{2}}U=L,$
say, is a unimodular lattice. Moreover since $\bar U$ is totally isotropic
then $L$ is integral and hence is a Niemeier lattice. Thus again we see
how to embed $U=\sqrt{2}L\to \L.$ Indeed, for 
a Niemeier lattice $L$ together with an isometric embedding $\sqrt{2}L\to
\L$ we have already observed that $|\L:\sqrt{2}L|=2^{12}$ and $\sqrt{2}L$ maps 
to a totally isotropic subspace of $\bar \L,$ hence is Lagrangian. We have
shown

\th Lemma 4.3

There is a natural bijection between embedded and re-scaled Niemeier lattices
$\sqrt{2}L\to \L$ and Lagrangian subspaces of $\L/2\L.$ 
\endth

Next we lift $\L/2\L$ to an extra-special group $Q\simeq 2_+^{1+24}.$
As is well-known [G], the subgroup $C$ of the Monster which leaves invariant
the ``untwisted'' part $V^+_{\L}$ of $V^{\natural}$ (cf. (4.1)) contains
$Q$ as a normal subgroup with quotient $C/Q\simeq Co_1.$

The Lagrangian subspaces of $\bar \L$ are precisely those which lift
to (maximal) elementary abelian subgroups of $Q$ isomorphic to
$\Z_2^{13}.$ Now in the Monster there are two classes of involutions,
of types $2A$ ($2+$) and $2B$ ($2-$) respectively, and one may ask how
they distribute themselves among the $\Z_2^{13}$ subgroups. There is a
pretty answer, which runs as follows: one knows (cf. [C1], for example)
that the $2A$ involutions of $Q$ map onto elements of $\L/2\L$ which 
themselves lift to elements $\l\in\L$ satisfying $(\l,\l)=4,$
whereas the $2B$ involutions correspond to $\l$ such that $(\l,\l)=8.$ Thus
for a fixed maximal elementary abelian subgroup $E\leq Q$ corresponding
to the re-scaled Niemeier $\sqrt{2}L\subset \L,$ the elements in $E$
of type $2A$ correspond precisely to the elements of squared-length 2, i.e.,
to the elements of $L_2.$ Of course, these form a Niemeier root system. So we
have proved

\th Proposition 4.4

Let $E\subset Q$ be a maximal elementary abelian subgroup
which corresponds (via the bijection of Lemma 4.3) to the
re-scaled Niemeier lattice  $\sqrt{2}L.$ Then there is  a natural bijection
$$\{ 2A\ {\rm involutions\ in}\ E\}\to {\rm root\ system\ of\ type\ }
L_2.$$
\endth

Let $E$ be as in Proposition 4.4, corresponding to $\sqrt{2}L.$ We come full
circle with the next result,
which identifies the sub \voa of $V^{\natural}$ fixed
pointwise by $E.$ 

\th Proposition 4.5

There is a natural isomorphism of \voas 
$$(V^{\natural})^E\simeq V_{\sqrt{2}L}^+.$$
\endth

\pf Proof

Consider $V_{\L},$ which is linearly isomorphic to $S(H_-)\otimes \C[\L]$
(cf. [FLM]) where $H=\C\otimes \Z,$ $H_-=H_{-1}\oplus H_{-2}\oplus\cdots$
with each $H_i\simeq H,$ and $\C[\L]$  is the group algebra  of $\L.$ Now
the center $Z(Q)$ of $Q$  acts on $V^{\natural}$  with fixed-point sub \voa naturally isomorphic  to $V^+_{\L},$ so we can study  the action
of $E/Z(Q)$  on $V^+_{\L}$ to prove the Proposition.

Now elements $\gamma$ of $\R\otimes_{\Z}\L$ act on $V_{\L}$ by fixing $S(H_-)$
identically and acting on basis vectors $e^{\a},$ $\a\in\L,$ via
$$\gamma\cdot e^{\a}=e^{2\pi i(\gamma,\a)}e^{\a}.$$
This induces an action of the Leech torus $\R^{24}/\L$
on $V_{\L}$ in which $E$ corresponds to $(\sqrt{2}L)^*/{\L}={1\over \sqrt{2}}
L/\L.$ 
Thus if $\a\in\L,$ $\gamma\cdot e^{\a}=e^{\a}$ for all $\gamma\in {1\over \sqrt{2}}L$ if, and only if, $\a\in\sqrt{2}L.$

This shows that $V_{\L}^{{1\over \sqrt{2}}L/\L}\simeq V_{\sqrt{2}L},$ and 
therefore also 
$$(V^{\natural})^E\simeq  (V_{\L}^+)^{{1\over \sqrt{2}}L/\L}\simeq (V_{\L}^{{1\over \sqrt{2}}L/\L})^+\simeq V^+_{\sqrt{2}L}.$$
\qed

Finally, let us continue to let $E$ and $L$ be as above. Restricting
Proposition 4.5 to the Griess algebra yields an isomorphism of algebras
$$B^E\simeq  (V^+_{\sqrt{2}L})_2$$
and in particular $B^E\supset B^+(L_2)=(V^+_{\sqrt{2}L})_2$ where  $R$ is the
root lattice $L_2.$  This refines the containment of
$B^+(L_2)$ in $B$ found earlier.

We have seen in Theorem 3.1 (and its extension to semi-simple root
systems) that there is a surjection of $A(L_2)$ onto $B^+(L_2)$ and
in particular $B^+(L_2)$ is generated by the images of
$t(\a),u(\a)$ for $\a\in L_2^+.$

Now as $t(\a)^2=8t(\a)$ and $\<t(\a),t(\a)\>=4$ then each 
$t(\a)/8$ corresponds to an idempotent with central charge $1/2$ (cf. (2.4)).
Of course $\a$ itself corresponds to an involution of type $2A.$ It
turns out that the bijection $\a\mapsto t(\a)/8$ is precisely the 
correspondence between {\it transpositions}
(i.e., $2A$ involutions of the Monster) and {\it transposition
axes} in the Griess algebra. See [C2] and [M] for more information
on this point. Using this perspective, one can read off the relations 
satisfied by the images of the $t(\a)/8,$ from Table 1 of [Nor].
This was the original motivation for introducing the algebra $A(\P).$

 \Bib        Bibliography
 
 
\rf {B} N. Bourbaki, Groupes et ${\rm alg}\grave{\rm e}{\rm bres}$ de Lie, Chaps. 4, 5, 6,
Hermann, Paris, 1968.

 \rf {C1}   J. H. Conway, Three lectures on exceptional groups, Chap. 10 of
J. H. Conway and N. J. A. Sloane,
{\it Sphere Packing, Lattices and Geoups,} Springer-Verlag, New York.

\rf {C2}    J. H. Conway, A simple construction for the Fischer-Griess monster group, {\it Invent. Math.} {\bf 79} (1985), 513-540.
 
\rf {Cur} R. T. Curtis, On subgroups of $\cdot 0$. II. Local Structure,
{\it J. Algebra,} {\bf 63} (1980), 413-434.

 \rf {D} C. Dong, Representations of the moonshine module 
vertex operator algebra, {\it Contemporary Math.} {\bf 175} (1994),
27-36.

  \rf {DL} C. Dong, J. Lepowsky, 
Generalized Vertex
Algebras and Relative Vertex Operators, {\it Progress in Math.} Vol. 112,
Birkh\"{a}user, Boston 1993.

\rf {DLM} C. Dong, H. Li, G. Mason, 
Some twisted modules for the moonshine vertex operator algebras,
{\it Contemp. Math.} {\bf 193} (1996), 25-43.

\rf {DMZ} C. Dong, G. Mason, Y. Zhu, 
Discrete series of the 
Virasoro algebra and the moonshine module, {\it Proc. Symp. Pure. Math., American Math. Soc.} {\bf 56} II (1994), 295-316.

\rf {FLM} I. B. Frenkel, J. Lepowsky, A. Meurman,
Vertex Operator Algebras and the Monster, {\it Pure and Applied
Math.,} Academic Press (1988).

\rf {FQS} D. Friedan, Z. Qiu and S. Shenker, Conformal invariance,
unitarity and two-dimensional critical exponents, MSRI publ. \#3, Springer-Verlag,  419-449 (1985).
 
\rf {GKO} P. Goddard, A. Kent and D. Olive,
Unitary representations of the Virasoro Algebra and super-Virasoro algebras,
{\it Commun. Math. Phys.} 103, 105-119 (1986).

\rf {G} R. Griess Jr., The Friendly Giant,
 {\it Invent. Math.} {\bf 69} (1982), 1-102.

\rf {H} Y. Huang, A non-meromorphic extension of the
moonshine module vertex operator algebra, {\it Contemporary Math.}
{\bf 163} (1996), 123-148.

\rf {MN} W. Meyer and W. Neutsch, Associative subalgebras of the
Griess algebra, {\it J. Algebra} {\bf 158} (1993),1-17.

\rf {M} M. Miyamoto, Griess algebras and conformal vectors 
in vertex operator algebras,
{\it J. Algebra,} {\bf 179} (1996), 523-548.

\rf {Nie} H. V. Niemeier, Definite quadratische Formen der Dimension 24
und Diskriminante 1, {\it J. Number Theory} {\bf 5} (1973), 142-178.

\rf {Nor} S. Norton, The monster algebra: some new formulae, {\it Contemp. 
Math.} {\bf 193} (1996), 297-306.

\rf {V} B. B. Venkov, The classification of integral even unimodular 24-dimensional quadratic forms,
{\it Proc. Steklov Inst. Math.} {\bf 4} (1980), 63-74.

\rf {ZF} A. B. Zamolodchikov and V. A. Fateev, 
Nonlocal (parafermion) currents in two-dimensional conformal quantum
field theory and self-dual critical points in $Z_N$-symmetric
statistical systems, {\it Sov. Phys., JETP} {\bf 62} (1985), 215-225.

\endBib

\Coordinates
 Department of Mathematics\\
 University of California\\
 Santa Cruz, CA 95064 (C.D., H.L., G.M.)

 Email: dong@cats.ucsc.edu, hli@cats.ucsc.edu, gem@cats.ucsc.edu
 \endCoordinates
 
\Coordinates
DPMMS\\
16 Mill Lane\\
Cambridge CB2 1SB\\
U.K. (S.N.)
  
 Email: simon@pmms.cam.ac.uk
 \endCoordinates
 \end